\documentclass{article}
\usepackage{times}

\usepackage{theorem}
\newtheorem{theorem}{Theorem}[section]

\usepackage{amsmath}
\usepackage{amssymb}

\begin{document}

\date{}
\title{Cubic-matrix splines and \\ second-order matrix models}
\author{M.M.~Tung,
        L.~Soler,
        E.~Defez, and
        A.~Herv\'as\footnote{
        \texttt{\{mtung, edefez, ahervas\}@imm.upv.es}} \\
        Instituto de Matem\'atica Multidisciplinar \\
        Universidad Polit\'ecnica de Valencia, Spain
          }
\maketitle

\begin{abstract}
\noindent
We discuss the direct use of cubic-matrix splines to obtain continuous
approximations to the unique solution of matrix models of the type
$Y''(x)=f(x, Y(x))$. For numerical illustration, an estimation of the
approximation error, an algorithm for its implementation, and an
example are given.
\end{abstract}

\section{Introduction.}\label{seccion1}
Matrix initial value problems of the form:
\begin{equation}\label{problem}
\left.  \begin{array}{rcl}
Y''(x) & = & f(x, Y(x))   \\
\\
Y(a)  & = & Y_{0} \ , \ Y'(a) \ = \ Y_{1}
\end{array}
 \right\}  \ a \leq x \leq b \ ,
\end{equation}
are frequently encountered in different fields of physics and
engineering (see {\it e.g.}\/ \cite{Zhang}). In the scalar case,
numerical methods for the calculation of approximate solutions of
(\ref{problem}) can be found in \cite{coleman}. For matrix problems,
linear multi-step matrix methods with constant steps have been
studied in \cite{15a}. Although in this case there exist {\it a
priori} error bounds for these methods (expressed as function of the
data problem), these error bounds are given in terms of an
exponential which depends on the integration step $h$. Therefore, in
practice, $h$ will take too small values. Problems of the type
(\ref{problem}) can be written as an extended first-order matrix
problem. Such a standard approach, however, involves an increase of the
computational cost caused by the increase of the problem
dimension. Recently, cubic-matrix splines were used
 in the resolution of first-order matrix differential
systems \cite{a5}, obtaining approximations that, among other
advantages, were of class $C^1$ in the interval $[a,b]$, and easy to
compute producing an approximation error $O(h^4)$. The present
work extends this powerful scheme to the solution of matrix
problems of type (\ref{problem}). Throughout this work, we will
adopt the notation for norms and matrix cubic splines as
in~\cite{a5} and common in matrix calculus. The paper is organized
as follows. Section \ref{seccion2} develops the proposed method.
Finally, in Section \ref{seccion4}, an example is presented.
\section{Construction of the method.}\label{seccion2}
Let us consider the initial value problem
\begin{equation}\label{tema2}
\left.  \begin{array}{rcl}
Y''(x) & = & f(x, Y(x))   \\
\\
Y(a)  & = & Y_{0} \ , \ Y'(a) \ = \ Y_{1}
\end{array}
 \right\}  \ a \leq x \leq b \ ,
\end{equation}
where $Y_{0}, Y_{1}, Y(t) \in {\mathbb C}^{r \times q}$, $f:[a,b]
\times {\mathbb C}^{r \times q}\times \longmapsto {\mathbb C}^{r
\times q}$, $f \in {\cal C}^{0}\left(T \right)$, with
\begin{equation}\label{guapa1}
T \ = \ \left\{(x,Y) \ ; \ a \leq x \leq b \ , \ Y \in {\mathbb
C}^{r \times q} \right\} \ ,
\end{equation}
and $f$ fulfills the global Lipschitz's condition
\begin{equation}\label{2}
\left\|f\left(x, Y_{1}\right) \ - \ f\left(x, Y_{2}\right)\right\| \
\leq L \left\|Y_{1} - Y_{2} \right\| \ , \ a \leq x \leq b \ ,
Y_{1}, Y_{2} \in {\mathbb C}^{r \times q} \ .
\end{equation}
Let us also use the partition of the interval $[a,b]$ defined by
\begin{equation}\label{partition}
\Delta_{[a,b]}=\left\{a=x_{0}<x_{1}<\ldots <x_{n}=b \right\} \ , \
x_{k}=a+ k h \ , \ k=0,1,\ldots,n \ ,
\end{equation}
where $h=(b-a)/n$, $n$ being a positive integer. We will construct
in each subinterval $[a+k h, a+(k+1)h]$ a matrix-cubic spline
approximating the solution of problem (\ref{tema2}). For the first
interval $[a,a+h] $, we consider that the matrix-cubic spline is given
by
\begin{equation}\label{tema5}
S_{ \left|_{\left[a,a+h \right]}\right.}  (x) \ = \ Y(a) \ + \
Y'(a)(x-a) \ + \ \frac{1}{2!} Y''(a) (x-a)^2 \ + \ \frac{1}{3!}
A_{0} (x-a)^3 \ ,
\end{equation}
where $A_{0} \in {\mathbb C}^{r \times q}$ is a matrix parameter to
be determined. It is straightforward to check:
$$
S_{ \left|_{\left[a,a+h \right]}\right.}(a)=Y(a) \ , \
S'_{\left|_{\left[a,a+h \right]}\right.}(a)=Y'(a) \ , \
S''_{\left|_{\left[a,a+h \right]}\right.}(a)=Y''(a)=f(a,S_{
\left|_{\left[a,a+h \right]}\right.}(a)) \ .
$$
Thus, (\ref{tema5}) satisfies the equations of problem (\ref{tema2})
at point $x=a$.
To fully construct the matrix-cubic spline, we must still determine
$A_{0}$. By imposing that (\ref{tema5}) is a solution of problem
(\ref{tema2}) in $x=a+h$, we have:
\begin{equation}\label{3}
S''_{ \left|_{\left[a,a+h \right]}\right.} (a+h) \ = \ f\left(a+h,
S_{ \left|_{\left[a,a+h \right]}\right.}  (a+h) \right) \ ,
\end{equation}
and obtain from (\ref{3}) the matrix equation with only one unknown
matrix $A_{0}$:
\begin{equation}\label{4}
A_{0} \ = \ \frac{1}{h}\left[f\left(a+h,Y(a)+ Y'(a)h + \frac{1}{2}
Y''(a) h^2 + \frac{1}{6} A_{0} h^3 \right) - Y''(a) \right] \ .
\end{equation}
Assuming that the matrix equation (\ref{4}) has only one solution
$A_{0}$, the matrix-cubic spline is totally determined in the
interval $[a,a+h]$. Now, in the next interval $[a+h, a+2h]$, the
matrix-cubic spline is defined by:
\begin{eqnarray}\label{tema10}
S_{ \left|_{\left[a+h, a+2h \right]}\right.}  (x) & = &
S_{ \left|_{\left[a, a+h \right]}\right.}  (a\!+\!h)\!+\!
S'_{ \left|_{\left[a, a+h \right]}\right.}  (a\!+\!h) (x-(a+h)) \nonumber \\
& + & \frac{1}{2!} S''_{ \left|_{\left[a, a+h \right]}\right.}
(a\!+\!h) (x\!-\!(a\!+\!h))^2\!+\!\frac{1}{3!} A_{1}
(x\!-\!(a\!+\!h))^3 \ ,
\end{eqnarray}
so that $S(x)$ is of class ${\cal C}^2([a,a+h] \cup [a+h,a+2h])$,
and all of the coefficients of matrix-cubic spline $S_
{\left|_{\left[a+h, a+2h \right]} \right.} (x)$ are determined with
the exception of $A_{1} \in {\mathbb C}^{r \times q}$. By
construction, matrix-cubic spline (\ref{tema10}) satisfies the
differential equation (\ref{tema2}) in $x=a+h$. We can obtain
$A_{1}$ by requiring that the differential equation (\ref{tema2})
holds at point $x=a+2h$:
$$
S''_{ \left|_{\left[a+h,a+2h \right]}\right.} (a+2h) \ = \
f\left(a+2h, S_{ \left|_{\left[a+h,a+2h \right]}\right.}  (a+2h)
\right) \ .
$$
Expanding, we obtain the matrix equation with only one unknown
matrix $A_{1}$:
\begin{eqnarray}\label{4a}
A_{1}&=&\frac{1}{h}\!\left[f\left(a\!+\!2h, S_{ \left|_{\left[a,a+h
\right]}\right.}  (a\!+\!h)\!+\!S'_{ \left|_{\left[a,a+h
\right]}\right.} (a\!+\!h) h\!+\!\frac{1}{2} S''_{
\left|_{\left[a,a+h \right]}\right.} (a\!+\!h)
h^2\!+\!\frac{1}{6} A_{1} h^3\right) \right.   \nonumber \\
&-& \left. S''_{ \left|_{\left[a,a+h \right]}\right.}
(a\!+\!h)\!\right] \ .
\end{eqnarray}
Let us assume that the matrix equation (\ref{4a}) has only one
solution $A_{1} $. This way the spline is now totally determined in the
interval $[a+h,a+2h] $. Iterating this process, let us construct the
matrix-cubic spline taking $\left[a+(k-1) h, a+k h \right]$ as the
last subinterval. For the next subinterval $\left[a+kh,
a+(k+1)h\right]$, we define the corresponding matrix-cubic spline as
$$
S_{ \left|_{\left[a+k h, a+(k+1)h \right]}\right.}(x)=
\beta_{k}(x)+\frac{1}{3!}A_{k}(x-(a+k h))^3
$$
\begin{equation}\label{tema14}
\mbox{where} \ \ \beta_{k}(x)= \sum_{l=0}^{2} \frac{1}{l!}S^{(l)}_{
\left|_{\left[a+(k-1)h, a+kh \right]}\right.}(a+k h)(x-(a+k h))^{l}
\ .
\end{equation}
With this definition, it is $S(x) \in {\cal C}^2 \displaystyle \left
(\bigcup_{j=0}^{k} [a+jh,a+(j+1)h] \right) $ which fulfills the
differential equation (\ref{tema2}) at point $x=a+kh $. As an
additional requirement, we assume that $S(x)$ satisfies the
differential equation (\ref{tema2}) at the point $x=a+(k+1)h$,
{\it i.e.}
$$
S''_{ \left|_{\left[a+kh,a+(k+1)h \right]}\right.} (a+(k+1)h) \ = \
f\left(a\!+\!(k\!+\!1)h, S_{ \left|_{ \left[a+kh,a+(k+1)h
\right]}\right.} (a\!+\!(k\!+\!1)h) \right) \ .
$$
Subsequent expansion of this equation with the unknown matrix $A_{k}$ yields
\begin{equation}\label{tema15a}
A_{k} = \frac{1}{h} \left[ f\left( a+(k+1)h, \beta_{k}(a+(k+1)h)+
\frac{1}{6} A_{k} h^3\right) - \beta_{k}''(a+(k+1)h) \right] \ .
\end{equation}
Note that this matrix equation (\ref{tema15a}) is analogous to
equations (\ref{4}) and (\ref{4a}), when $k=0$ and $k=1$,
respectively. For a fixed $h$, we will consider the matrix function
of matrix variable $g:{\mathbb C}^{r \times q} \mapsto {\mathbb
C}^{r \times q}$ defined by
$$ 
g(T) = \frac{1}{h} \left[ f\left(
a+(k+1)h, \beta_{k}(a+(k+1)h)+ \frac{1}{6} T h^3\right) -
 \beta_{k}''(a+(k+1)h)\right] \ .
$$ 
Relation (\ref{tema15a}) holds if and only if
$A_{k}=g(A_{k})$, that is, if $A_{k}$ is a fixed point for function
$g(T)$. Applying the global Lipschitz's conditions (\ref{2}), it
follows that
$$
\left\|g(T_{1}) - g(T_{2}) \right\| \leq  \frac{L h^2}{6}
\left\|T_{1} - T_{2} \right\| \ .
$$
Taking $h<\sqrt{\frac{6}{L}}$, $g(T)$ yields a contractive matrix
function, which guarantees that equation (\ref{tema15a}) has unique
solutions $A_{k}$ for $k=0,1,\ldots,n-1$. Hence, the matrix-cubic
spline is now fully determined. Taking into account \cite[Theorem
5]{20}, the following result has been established:
\begin{theorem}\label{theorem}
If $h<\sqrt{\frac{6}{L}}$, then the matrix-cubic spline $S(x)$
exists in each subinterval $\left[a+kh, a+(k+1)h \right]$,
$k=0,1,\ldots,n-1$, as defined by the previous construction.
Furthermore, if $f\in {\cal C}^{1}(T)$, then $ \left\|Y(x) -
S(x)\right\| = O(h^3) \ \forall x \in [a,b],$ where $Y(x) $ is the
theoretical solution of system (\ref{tema2}).
\end{theorem}
Depending on the function $f$, matrix equations (\ref{4}) and
(\ref{tema15a}) can be solved explicitly or by using some iterative
method \cite{ort}. Summarizing , we have the following
\textbf{algorithm}:\newline

\begin{tabular}{l}
$\bullet$ Take $n>\displaystyle\frac{(b-a) \sqrt{L}}{\sqrt{6}}$, $h=(b-a)/n$ and 
 $\Delta_{[a,\ b]}$ defined by~(\ref{partition}). \\[.3cm]
$\bullet$ Solve (\ref{4}) and determine $S_ {\left|_{\left[a,a+h
\right]} \right.}(x)$
defined by~(\ref{tema5}). \\[.2cm]
$\bullet$ For $k=1$ to $n-1 $, solve (\ref{tema15a}). Determine $S_
{\left|_{\left[a+k h, a+(k+1)h \right]} \right.} (x) $
 defined by~(\ref{tema14}).
\end{tabular}
\section{Example}\label{seccion4}
The problem
\begin{equation}\label{ejemplo1}
Y''(t) \ + \ A Y(t)=0 \ ,
\end{equation}
with $Y(0)=Y_0$, $Y'(0)=Y_1$, has the exact solution
$$Y(t)=\cos{(\sqrt{A}t)} Y_0 +
\left(\sqrt{A}\right)^{-1}\sin{(\sqrt{A}t)} Y_1\ ,$$
where $\sqrt{A}$
denotes any square root of a non-singular matrix $A$,
\cite{hargreaves}. The principal drawback of this formal solution is the
difficult computation of $\sqrt{A}$, $\cos{(\sqrt{A}t)}$ and
$\sin{(\sqrt{A}t)}$. The proposed method avoids this drawback. We
consider problem (\ref{ejemplo1}) where $ A=\left(\begin{array}{cc}
1 &
0 \\
2 & 1
\end{array} \right), Y_0=\left(\begin{array}{cc} 0 & 0 \\
0 & 0 \end{array} \right) , Y_1=\left(\begin{array}{cc} 1 & 0 \\
1 & 1 \end{array} \right)$, $t\in [0,1]$, whose exact solution is
$Y(t)=\sin{\left[\left(\begin{array}{cc}
1 & 0 \\
1 & 1
\end{array} \right) t\right]}=\left(\begin{array}{cc}
\sin{(t)} & 0 \\
t \cos{(t)} & \sin{(t)}
\end{array} \right)$. In this case $ L \approx 2.82843 $.
By Theorem~\ref{theorem}, we
need to take $h<1.45647$, so we choose $h=0.1$ for example. The results
are summarized in the following table, where the numerical estimates
have been rounded to the fourth relevant digit. In each subinterval,
we evaluated the difference between the estimates of our numerical
approach and the exact solution. The maximum of these errors are
indicated in the third column.
\begin{center}
\begin{tabular}{||c|c|c||}
\hline
\mbox{Interval} & \mbox{Approximation} & \mbox{Max.\ Error} \\
\hline {\small $[0,0.1]$} & {\tiny $\left( \begin{array}{cc} x
-0.1664 x^3 & 0 \\   x - 0.4986 x^3 & x - 0.1664 x^3
\end{array}
\right)$ } & {\tiny $1.0072 \times 10^{-6}$}  \\
\hline {\small $[0.1,0.2]$} & {\tiny $\left(
\begin{array}{cc} 1.00005 x -0.0005 x^2 -
    0.1647 x^3 & 0 \\  1.0002 x - 0.0025 x^2 -
    0.4903 x^3 & 1.0001 x
-0.0005 x^2 -
    0.1647 x^3
\end{array}
\right)$ } & {\tiny $6.3032 \times 10^{-6}$}  \\
\hline {\small $[0.2,0.3]$} & {\tiny $\left(
\begin{array}{cc} 1.0005 x - 0.0025 x^2 -
    0.1614 x^3 & 0 \\  -0.0001 + 1.0022 x -
    0.0124 x^2 - 0.4738 x^3 & 1.0005 x - 0.0025 x^2 -
    0.1614 x^3
\end{array}
\right)$ } & {\tiny $2.0059 \times 10^{-5}$}  \\
\hline {\small $[0.3,0.4]$} & {\tiny $\left(
\begin{array}{cc} -0.0002 + 1.0018 x -
    0.0069 x^2 - 0.1565 x^3 & 0 \\ -0.0008 + 1.0088 x -
    0.0344 x^2 - 0.4494 x^3 & -0.0002 + 1.0018 x -
    0.0069 x^2 - 0.1565 x^3
\end{array}
\right)$ } & {\tiny $4.6213 \times 10^{-5}$}  \\
\hline {\small $[0.4,0.5]$} & {\tiny $\left(
\begin{array}{cc} -0.0006 + 1.0049 x -
    0.0147 x^2 - 0.1500 x^3 & 0 \\ -0.0028 + 1.0242 x -
    0.0728 x^2 - 0.4174 x^3 & -0.0006 + 1.0049 x -
    0.0147 x^2 - 0.1500 x^3
\end{array}
\right)$ } & {\tiny $8.8359 \times 10^{-5}$}  \\
\hline {\small $[0.5,0.6]$} & {\tiny $\left(
\begin{array}{cc} -0.0016 + 1.0109 x -
    0.0266 x^2 - 0.1420 x^3 & 0 \\ -0.0077 + 1.0536 x -
    0.1316 x^2 - 0.3782 x^3 & -0.0016 + 1.0109 x -
    0.0266 x^2 - 0.1420 x^3
\end{array}
\right)$ } & {\tiny $1.4964 \times 10^{-4}$}  \\
\hline {\small $[0.6,0.7]$} & {\tiny $\left(
\begin{array}{cc} -0.0036 + 1.0210 x -
    0.0436 x^2 - 0.1327 x^3 & 0 \\ -0.0176 + 1.1030 x -
    0.2140 x^2 - 0.3324 x^3 & -0.0036 + 1.0210 x -
    0.0436 x^2 - 0.1327 x^3
\end{array}
\right)$ } & {\tiny $2.3267 \times 10^{-4}$}  \\
\hline {\small $[0.7,0.8]$} & {\tiny $\left(
\begin{array}{cc} -0.0073 + 1.0368 x -
    0.0661 x^2 - 0.1219 x^3 & 0 \\ -0.0354 + 1.1791 x -
    0.3227 x^2 - 0.2807 x^3 & -0.0073 + 1.0368 x -
    0.0661 x^2 - 0.1219 x^3
\end{array}
\right)$ } & {\tiny $3.3941 \times 10^{-4}$}  \\
\hline {\small $[0.8,0.9]$} & {\tiny $\left(
\begin{array}{cc} -0.0134 + 1.0597 x -
    0.0947 x^2 - 0.1100 x^3 & 0 \\ -0.0646 + 1.2885 x -
    0.4595 x^2 - 0.2237 x^3 & -0.0134 + 1.0597 x -
    0.0947 x^2 - 0.1100 x^3
\end{array}
\right)$ } & {\tiny $4.7114 \times 10^{-4}$}  \\
\hline {\small $[0.9,1]$} & {\tiny $\left(
\begin{array}{cc} -0.0229 + 1.0914 x -
    0.1299 x^2 - 0.0970 x^3 & 0 \\ -0.1093 + 1.4378 x -
    0.6253 x^2 - 0.1623 x^3 & -0.0229 + 1.0914 x -
    0.1299 x^2 - 0.0970 x^3
\end{array}
\right)$ } & {\tiny $6.2838 \times 10^{-4}$}  \\
\hline \hline
\end{tabular}
\end{center}

\end{document}